\documentclass[12pt]{article}
\usepackage{amsmath,graphicx,amssymb,amsthm}
\allowdisplaybreaks
\def\otime{\bar{\otimes}}

\def\AA{\mathfrak A}
\def\BB{\mathfrak B}

\def\A{\mathcal A}

\def\K{\mathcal K}

\def\B{\mathcal B}

\def\H{\mathcal H}
\def\P{\mathcal P}
\def\R{\mathcal R}
\def\S{\mathcal S}
\def\M{\mathcal M}
\def\N{\mathcal N}
\def\L{\mathcal L}
\def\Z{\mathcal Z}
\def\S{\mathcal S}

\def\amslatex{$\mathcal{A}\kern-.1667em\lower.5ex\hbox{$\mathcal{M}$}\kern-.125em\mathcal{S}$-\LaTeX}

\newtheorem{set}{set}[section]
\newtheorem{Corollary}[set]{Corollary}

\newtheorem{Lemma}[set]{Lemma}

\newtheorem{Theorem}[set]{Theorem}
\newcommand{\define}{\mathrel{\hbox{$\equiv$\hskip -.90em \lower .47ex \hbox{$\leftharpoondown$}}}}
\newcommand{\enifed}{\mathrel{\hbox{$\equiv$\hskip -.90em \lower .47ex \hbox{$\rightharpoondown$}}}}

\begin{document}
\pagestyle{myheadings}
\date{}
\title{On maximal injective subalgebras of tensor products of von
Neumann algebras}
\author{Junsheng Fang\sc \footnote{jfang@cisunix.unh.edu.}\\
\emph{\small Department of Mathematics, University of New Hampshire,
Durham, NH 03824, USA}}
 \maketitle

\begin{abstract}Let $\M_i$ be a
von Neumann algebra, and $\BB_i$ be a maximal injective von Neumann
subalgebra of $\M_i$, $i=1,2$. If $\M_1$ has separable predual and
the center of $\BB_1$ is atomic, e.g., $\BB_1$ is a factor, then
$\BB_1\otime\BB_2$ is a maximal injective von Neumann subalgebra of
$\M_1\otime\M_2$. This partly answers a question of Popa.
\end{abstract}

\vskip1.0cm
\centerline{{\Large\bf Introduction}} \mbox{}\\

F.J.~Murray and J.von~Neumann~\cite{M-V1,M-V2,M-V3,vN1,vN2}
introduced and studied certain algebras of Hilbert space operators.
Those algebras are now called ``Von Neumann algebras." They are
strong-operator closed self-adjoint subalgebras of all bounded
linear transformations on a Hilbert space. \emph{Factors} are von
Neumann algebras whose centers consist of scalar multiples of the
identity. Every von Neumann algebra is a direct sum (or ``direct
integral'') of factors. Thus factors are the building blocks for all
von Neumann algebras.

Murray and von Neumann~\cite{M-V1} classified factors by means of a
relative dimension function. \emph{Finite factors} are those for
which this dimension function has range the closed interval $[0,c]$
for some positive $c$. For finite factors, this dimension function
gives rise to a (unique, when normalized) tracial state. In general,
a von Neumann algebra admitting a faithful normal trace is said to
be \emph{finite}. Finite-dimensional ``finite factors" are full
matrix algebras $M_n(\mathbb{C})$, $n=1,2,\cdots$.
Infinite-dimensional ``finite factors" are called factors of type
$II_1$. \emph{Infinite factors} are those for which the range of the
dimension function includes $\infty$.

In~\cite{M-V3}, Murray and von Neumann introduced and studied a
family of factors of type $II_1$ very closely related to matrix
algebras. Murray and von Neumann called these factors
\emph{approximately finite} since they are the ultraweak closure of
the ascending union of a family of finite-dimensional self-adjoint
subalgebras. They proved that all ``approximately finite" factors of
type $II_1$ are * isomorphic.
 Since these factors are finite, J.~Dixmier~\cite{Di} considered the term ``approximately
finite" inappropriate and called them \emph{hyperfinite}.  However,
for infinite factors possessing the same property, the term
``hyperfinite" is also inappropriate. So later on the  name
\emph{approximately finite dimensional} (AFD) were introduced for
these factors.

A von Neumann algebra $\BB$ acting on a Hilbert space $\H$ is called
\emph{injective} if there is a norm one projection from $\B(\H)$,
the algebra of all bounded linear operators on $\H$, onto $\BB$.
Since the intersection of a decreasing sequence of injective
algebras is injective, and the commutant of an injective algebra is
injective, every AFD factor is injective. In~\cite{Co}, A.~Connes
proved that  a separable injective von Neumann algebra (von Neumann
algebra with separable predual) is approximately finite dimensional.
As a corollary, this shows that the hyperfinite type $II_1$ factor
$\R$ is the unique separable injective factor of type $II_1$. The
proof of connes' result is so deep and rich in ideas and techniques
that it remains a basic resource in the subject.

Compare with injective factors, non-injective factors (even
non-injective type $II_1$ factors) are far from being understood. A
standard method of investigation in the study of general factors is
to study the injective von Neumann subalgebras of these factors.
Along this line, we have R.~Kadison's question (Problem 7
in~\cite{Ka}): Does each self-adjoint operator in a $II_1$ factor
lie in some hyperfinite subfactor? Since every separable abelian von
Neumann algebra is generated by a single self-adjoint operator,
Kadison's question has an equivalent form: Is each separable abelian
von Neumann algebra of a $II_1$ factor contained in some hyperfinite
subfactor?

Let $\M$ be a type $II_1$ factor with a faithful normal trace
$\tau$. If $T=T^*=\sum_{k=1}^n\lambda_kE_k$ is a self-adjoint
operator in $\M$  such that $\sum_{k=1}^nE_k=I$ and $\tau(E_k)=1/n$
for $1\leq k\leq n$, then $T$ is in a type $I_n$ subfactor $\M_n$ of
$\M$ which has $E_1,E_2,\cdots,E_n$ as diagonals. Since the set
$\S=\{T:\, T=T^*=\sum_{k=1}^n\lambda_kE_k\,\, \text{such that} \,\,
\sum_{k=1}^nE_k=I \,\, \text{and}\,\,
\tau(E_1)=\cdots=\tau(E_n)=\frac{1}{n},\, n=1,2,\cdots\}$ is dense
in the set of self-adjoint operators in $\M$ relative to the
strong-operator topology, for each self-adjoint operator $T$ in $\M$
we can choose a sequence $\{T_n\}_{n=1}^{\infty}\subseteq \S$ such
that $T_n$ converges to $T$ in the strong operator topology. So one
may expect that the answer to Kadison's question could be
affirmative if one very carefully constructs $T_n$ and $\M_n$ for
each $n$.

Out of expectation, this problem was answered in the negative in a
remarkable paper~\cite{Po2} by S.~Popa. In~\cite{Po2}, Popa showed
that if $\L(\mathbb{F}_n)$ is the type $II_1$ factor associated with
the left regular representation $\lambda$ of the free group
$\mathbb{F}_n$ on $n$ generators, $2\leq n\leq \infty$, and $a$ is
one of the generators of $\mathbb{F}_n$, then the abelian von
Neumann subalgebra generated  by the unitary $\lambda(a)$ is a
maximal injective von Neumann subalgebra of $\L(\mathbb{F}_n)$. So
quite surprisingly, a diffuse abelian von Neumann algebra can be
embedded in a type $II_1$ factor as a maximal injective von Neumann
subalgebra!

By considering actions of free groups on non-atomic probability
spaces, Popa constructed  more examples of maximal injective von
Neumann subalgebras in factors of type $II_1$. In~\cite{Ge}, L.~Ge
showed that every non-atomic injective von Neumann algebra with
separable predual is maximal injective in its free product with any
von Neumann algebra associated with a countable discrete group. Popa
raised the following question in~\cite{Po2}: If $\M_1,\M_2$ are type
$II_1$ factors and $\BB_1\subseteq\M_1,\BB_2\subseteq\M_2$ are
maximal injective von Neumann subalgebras, is $\BB_1\otime\BB_2$
maximal injective in $\M_1\otime\M_2$? He also asked if this is true
when we assume that $\BB_1=\M_1$ is the hyperfinite $II_1$ factor.
This question, if answered in the affirmative, would considerably
enlarge our class of examples.

 The first break of Popa's question was obtained by Ge and Kadison~\cite{Ge,G-K}. By applying
their remarkable ``splitting theorem", Ge and Kadison answered the
second part of Popa's question affirmatively. Precisely, Ge and
Kadison proved that if $\M_1$ is an injective factor, and  $\BB_2$
is a maximal injective von Neumann subalgebra of $\M_2$, then
$\M_1\otime\BB_2$ is maximal injective in $\M_1\otime\M_2$.
In~\cite{S-Z}, S.~Str\v{a}til\v{a} and L.~Zsid\'{o} improved
Ge-Kadison's result by removing the factor condition of $\M_1$: if
$\M_1$ is an injective von Neumann algebra and $\M_2$ is a von
Neumann algebra with separable predual, and $\BB_2$ is a maximal
injective von Neumann subalgebra of $\M_2$, then $\M_1\otime\BB_2$
is maximal injective in $\M_1\otime\M_2$.

In this paper we answer Popa's question affirmatively in a more
general setting. Our main result is the following. Let $\M_i$ be a
von Neumann algebra, and $\BB_i$ be a maximal injective von Neumann
subalgebra of $\M_i$, $i=1,2$. If $\M_1$ has separable predual and
the center of $\BB_1$ is atomic, e.g., $\BB_1$ is a factor, then
$\BB_1\otime\BB_2$ is maximal injective in $\M_1\otime\M_2$.

The paper is divided into six sections. Section 1 contains some
preliminaries and one key lemma. Another key lemma is proved in
section 2. Some direct applications are also given. In section 3, we
prove our main result in the special case when $\BB_1$ is a factor.
The main result  is proved in section 4. In section 5, we consider
the question: If $\M_1,\M_2$ are factors,
 and $\R_1,\R_2$ are maximal injective subfactors of $\M_1,\M_2$,
 respectively, is $\R_1\otime\R_2$ a maximal injective subfactor
 of $\M_1\otime \M_2$? We prove the following result: Let $\M_1,\M_2$ be factors,
 and $\R_1,\R_2$ be maximal injective subfactors of $\M_1,\M_2$,
 respectively. If $\R_1^{'}\cap\M_1\simeq\mathbb{C}^N$ ($1\leq N\leq
 \infty$)
  and $\R_2^{'}\cap\M_2=\mathbb{C}I$,
then $\R_1\otime\R_2$ is a maximal injective subfactor of
$\M_1\otime \M_2.$ In the last section, we mention some questions
related to Popa's question.

For the general theory of von Neumann algebras, we refer
to~\cite{Di,K-R}, and~\cite{St}.

The author wishes to express his deep gratitude to professor Liming
Ge and professor Don Hadwin, for valuable conversations and their
encouragement.

\section{Preliminaries}
There are five topics in this section: injective von Neumann
algebras, maximal injective von Neumann subalgebras, minimal
injective von Neumann algebra extensions, maximal injective
subfactors, and two basic theorems on tensor products of von
Neumann algebras: Ge-Kadison's splitting theorem and Tomiyama's
slice mapping theorem. Lemma 1.2 is one of two key lemmas.
\subsection{Injective von Neumann algebras}
A \emph{conditional expectation} ${\bf E}$  from a von Neumann
algebra $\M$ onto a von Neumann subalgebra $\N$ is a positive,
linear mapping such that ${\bf E}(S_1TS_2)=S_1{\bf E}(T)S_2$ for all
$S_1,S_2$ in $\N$ and all $T$ in $\M$. J.~Tomiyama~\cite{To1} showed
that an idempotent of norm 1 from $\M$ onto $\N$ is a conditional
expectation. A von Neumann algebra
 $\BB$ acting on a Hilbert space $\H$ is
called injective if there is a conditional expectation from $\B(\H)$
onto $\BB$. If $\BB$ is  a von Neumann subalgebra of a von Neumann
algebra $\M$ and $\BB$ is injective, there is a conditional
expectation  from $\M$ onto $\BB$.

Let $\BB$ be a von Neumann algebra acting on a Hilbert space $\H$.
Then $\BB$ is injective if and only if the commutant $\BB^{'}$ of
$\BB$ is injective. Recall that if $E$ is a projection in $\BB$,
then the reduced von Neumann algebra of $\BB$ with respect to $E$
is the algebra $\BB_E\triangleq E\BB E.$ If $\BB$ is injective and
${\bf E}$ is a conditional expectation from $\B(\H)$ onto $\BB$,
then ${\bf E}$ induces a conditional expectation ${\bf E}_E$ from
$\B(E\H)$ onto $\BB_E$ by ${\bf E}_E(T)={\bf E}(ETE)$ for any
$T\in\B(E\H)$. Thus $\BB_E$ is an injective von Neumann algebra.

\subsection{Maximal injective von Neumann subalgebras}
Let $\M$ be a von Neumann algebra. A von Neumann subalgebra $\BB$
of $\M$ is called \emph{maximal injective} if it is injective and
if it is maximal with respect to inclusion in the set of all
injective von Neumann subalgebras of $\M$. If  $\{\BB_\alpha\}$ is
a family of  injective von Neumann subalgebras of $\M$ which is
inductively ordered by inclusion, then the weak operator closure
of $\bigcup_\alpha\BB_\alpha$ is an injective von Neumann
subalgebra of $\M$ which contains all $\BB_\alpha$. By Zorn's
lemma, $\M$ has maximal injective von Neumann subalgebras.

If $\M$ is a separable type $II_1$ factor, then $\M$ contains a
hyperfinite subfactor $\R$ such that $\R^{'}\cap\M=\mathbb{C}I$
(Corollary 4.1 of~\cite{Po1}). If $\BB$ is a maximal injective von
Neumann subalgebra of $\M$ which contains $\R$, then
$\BB^{'}\cap\M\subseteq\R^{'}\cap\M=\mathbb{C}I.$  In particular,
$\BB$ is an injective factor. By~\cite{Co}, $\BB$ is hyperfinite. So
every separable type $II_1$ factor  contains a hyperfinite subfactor
as a maximal injective von Neumann subalgebra.

In~\cite{Po2}, Popa exhibited  concrete examples of maximal
injective von Neumann subalgebras of type $II_1$ factors. Popa
showed that if $\L(\mathbb{F}_n)$ is the type $II_1$ factor
associated with the left regular representation $\lambda$ of the
free group $\mathbb{F}_n$ on $n$ ($2\leq n\leq \infty$) generators,
and if $a$ is one of the generators of $\mathbb{F}_n$, then the
abelian von Neumann algebra generated by the unitary $\lambda(a)$ is
a maximal injective von Neumann subalgebra of $\L(\mathbb{F}_n)$.
 In~\cite{Ge}, Ge showed that each non-atomic
injective von Neumann algebra with  separable predual is maximal
injective in its free product with any von Neumann algebra
associated with a countable discrete group. Note that any maximal
injective von Neumann subalgebra of a type $II_1$ factor must be
non-atomic.

If $\BB$ is a maximal injective von Neumann subalgebra of $\M$, then
$\BB$ is singular in $\M$, i.e., its normalizers in $\M$ are unitary
elements in $\BB$. Indeed, if $U$ is a unitary element in $\M$ and
$U\BB U^*=\BB$, then the von Neumann subalgebra of $\M$ generated by
$\BB$ and $U$ is also injective. Since $\BB$ is maximal injective in
$\M$, $U\in\BB$. In particular, it follows that
$\BB^{'}\cap\M\subseteq \BB$. Let $\Z$ be the center of $\BB$. We
have $\Z\subseteq\BB^{'}\cap\M\subseteq\BB^{'}\cap\BB=\Z$, which
implies that $\Z=\BB^{'}\cap\M$. We summarize these facts in the
following lemma.
\begin{Lemma} Let $\BB$ be a maximal injective von Neumann
subalgebra of $\M$. Then $\BB$ is singular in $\M$. In particular,
$\Z=\BB^{'}\cap\BB=\BB^{'}\cap\M$.
\end{Lemma}

\subsection{Minimal injective von Neumann algebra extensions}

Let $\N$ be a von Neumann algebra. An injective von Neumann
algebra $\AA$ is called a \emph{minimal injective von Neumann
algebra extension} of $\N$ if $\AA\supseteq\N$ and if  it is
minimal with respect to inclusion in the set of all injective von
Neumann algebras which contain $\N$.

Let $\M$ be a von Neumann algebra acting on a Hilbert space $\H$ and
$\BB$ be a maximal injective von Neumann subalgebra of $\M$. Then
$\BB^{'}$, the commutant of $\BB$, is a minimal injective von
Neumann algebra extension of $\M^{'}$. Indeed, if $\L^{'}$ is an
injective von Neumann algebra such that
$\M^{'}\subseteq\L^{'}\subseteq\BB^{'}$, then $\L=(\L^{'})^{'}$ is
an injective von Neumann algebra  such that
$\BB\subseteq\L\subseteq\M$. Since $\BB$ is a maximal injective von
Neumann subalgebra of $\M$, $\BB=\L$. By von Neumann's double
commutant theorem~\cite{vN1}, $\BB^{'}=\L^{'}$.

Let  $\AA$ be a minimal injective von Neumann algebra extension of a
von Neumann algebra $\N$. Let  $\varphi$ be a faithful normal
representation of $\AA$ on a Hilbert space $\H$. Then
$(\varphi(\AA))^{'}$ is a maximal injective von Neumann algebra of
$(\varphi(\N))^{'}$. Indeed, if $\L$ is an injective von Neumann
algebra such that
$(\varphi(\AA))^{'}\subseteq\L^{'}\subseteq(\varphi(\N))^{'}$, then
$\varphi(\N)\subseteq\L\subseteq\varphi(\AA)$ and $\L$ is injective.
Thus $\N\subseteq\varphi^{-1}(\L)\subseteq\AA$ and
$\varphi^{-1}(\L)$ is injective. Since $\AA$ is a minimal injective
von Neumann algebra extension of $\N$, $\varphi^{-1}(\L)=\AA$. Hence
$\L=\varphi(\AA)$ and
$\L^{'}=(\varphi(\AA))^{'}$.\\

In~\cite{Ha}, U.~Haagerup proved that any von Neumann algebra is
*-isomorphic to a von Neumann algebra $\M$ on a Hilbert space
$\H$, such that there is a conjugate linear, isometric involution
$J$ of $\H$ and a self-dual cone $\P$ in $\H$ with the properties:
\begin{enumerate}
\item $J\M J=\M^{'}$, \item $JZJ=Z^*$, for $Z$ in the center of
$\M$, \item $J\xi=\xi,\,\, \xi\in\P$, \item
$XJXJ(\P)\subseteq\P\,\,\text{for all}\,\, X\in\M$.
\end{enumerate}
A quadruple $(\M,\H,J,\P)$ satisfying the conditions 1-4 is called
a \emph{standard form} of the von Neumann algebra $\M$.
 Recall that a von Neumann algebra $\M$ acting on a Hilbert space $\H$ is said
to be \emph{standard} if there exists a conjugation
$J:\H\rightarrow\H$, such that the mapping $X\rightarrow JX^*J$ is a
* anti-isomorphism from $\M$ onto $\M^{'}$. If $\M$ is standard on
$\H$, we can choose $J$ and $\P$ in $\H$, such that $(\M,\H,J,\P)$
is a standard form (cf~\cite{Ha}, Theorem 1.1).  Let $\M$ be
standard on $\H$, and $\theta$ be a *-automorphism of $\M$, then
there is a unitary operator $U$ on $\H$ such that $\theta(X)=UXU^*$
for all $X\in\M$ (cf~\cite{Ha}, Theorem 3.2).

 The following
lemma, which has an independent interest, is a key lemma.

\begin{Lemma}Let $\N$ be a von Neumann algebra and  $\AA$ be a
minimal injective von Neumann algebra extension of $\N$. If
$\theta\in Aut(\AA)$ \emph{(}the group of all *-automorphisms of
$\AA$\emph{)} satisfies $\theta(X)=X$ for all $X\in\N$. Then
$\theta(Y)=Y$ for all $Y\in\AA$.
\end{Lemma}
\begin{proof} We can assume that $\AA$ is standard on a Hilbert
space $\H$. Then there is a unitary operator $U\in\B(\H)$ such that
$\theta(Y)=UYU^*$ for all $Y\in\AA$. Since for all $X\in\N$ we have
$\theta(X)=X$, $UXU^*=X$. Thus, $U\in\N{'}$. Define
$\theta'(Y')=UY'U^*$ for $Y'\in\AA^{'}$. Note that for all $Y\in\AA,
Y'\in\AA^{'}$, $\theta'(Y')\theta(Y)=\theta(Y)\theta'(Y')$. Since
$\theta(\AA)=\AA$, $\theta'(Y')\in\AA^{'}$ for all $Y'\in\AA^{'}$,
i.e., $U\AA^{'} U^*\subseteq\AA^{'}$. Note that
$\theta^{-1}(Y)=U^*YU$ is also a *-isomorphism of $\AA$. Same
arguments as above show that $U^*\AA^{'} U\subseteq \AA^{'}$. So
$U\AA^{'} U^*=\AA^{'}$. This implies that $U\in\N{'}$ is in the
normalizer of $\AA^{'}$. Note that $\AA^{'}$ is maximal injective in
$\N{'}$. By Lemma 1.1, $U\in\AA^{'}$. So $\theta(Y)=UYU^*=Y$ for all
$Y\in\AA$.
\end{proof}
\begin{Lemma} If $\AA$ is a minimal injective von
Neumann algebra extension of a von Neumann algebra $\N$, and $P,Q$
are non-zero central projections in $\AA$ such that $PQ=0$, then
there does not exist a *-isomorphism $\phi$ from $\AA_P$ onto
$\AA_Q$ such that $\phi(PX)=QX$ for all $X\in\N$.
\end{Lemma}
\begin{proof} Otherwise, assume $\phi$ is a *-isomorphism from $\AA_P$
onto $\AA_Q$ such that $\phi(PX)=QX$ for all $X\in\N$.  For any
$Y\in\AA$, $Y=PY+QY+(I-P-Q)Y$. Define $\theta$ from $\AA$ to $\AA$
by: $\theta(Y)=\phi(PY)+\phi^{-1}(QY)+(I-P-Q)Y$. Since $P,Q$ are
mutually orthogonal central projections in $\AA$ and $\phi$ is a
*-isomorphism from $\AA_P$ onto $\AA_Q$, $\theta\in Aut(\AA)$.
Note that for any $X\in\N$,
$\theta(X)=\phi(PX)+\phi^{-1}(QX)+(I-P-Q)X=QX+PX+(I-P-Q)X=X$. Since
$\AA$ is a minimal injective von Neumann algebra extension of $\N$,
by Lemma 1.2, $\theta(Y)=Y$ for all $Y\in\AA$. Therefore,
$P=\theta(P)=\phi(P)=Q$. Now we have $P=PQ=0$. It contradicts to the
assumption that $P\neq 0$.
\end{proof}
\begin{Corollary}Let $\AA$ be a minimal injective von
Neumann algebra extension of a von Neumann algebra $\N$, and $P,Q$
be  central projections in $\AA$. If there is a *-isomorphism $\phi$
from $\AA_P$ onto $\AA_Q$ such that $\phi(PX)=QX$ for all $X\in\N$,
then $P=Q$ and $\phi(PY)=PY$ for all $Y\in\AA$.
\end{Corollary}
\begin{proof} Suppose $P\neq Q$. Let $R=PQ$. Without loss of generality, we can
assume that $P_1=P-R>0$. Let $Q_1=\theta(P_1)\leq Q$.  Then $P_1,
Q_1$ are non-zero central projections in $\AA$ and
$P_1Q_1=Q_1P_1=0$. Since $\phi$ is a *-isomorphism from $\AA_P$ onto
$\AA_Q$, $\phi$ induces a *-isomorphism $\psi$ from $\AA_{P_1}$ onto
$\AA_{Q_1}$ such that $\psi(P_1Y)=\phi(P_1Y)$ for all $Y\in\AA$.
Since for any $X\in\N$, $\psi(P_1X)=\phi(P_1X)=\phi(P_1)\phi(PX)=Q_1
QX=Q_1X$. It contradicts to Lemma 1.3. Thus $P=Q$. Define
$\theta(Y)=\phi(PY)+(I-P)Y$, then $\theta\in Aut(\AA)$ and
$\theta(X)=X$ for any $X\in\N$. Since $\AA$ is a minimal injective
von Neumann algebra extension of $\N$, by Lemma 1.2, $\theta(Y)=Y$
for all $Y\in\AA$. Hence, $PY=\theta(PY)=\phi(PY)$ for all
$Y\in\AA$.
\end{proof}
\subsection{Maximal injective subfactors}

In~\cite{F-K},B.~Fuglede and Kadison established the existence of
\emph{maximal hyperfinite subfactors} of a type $II_1$ factor. Since
a separable type $II_1$ factor is injective if and only if it is
hyperfinite, a subfactor of a separable type $II_1$ factor is a
maximal injective subfactor if and only if it is a maximal
hyperfinite subfactor. So every separable type $II_1$ factor has
maximal injective subfactors.

In~\cite{F-K}, Fuglede and Kadison also asked if each maximal
hyperfinite subfactor of a $II_1$ factor has a trivial relative
commutant (that is, only the scalars in the factor commute with the
subfactor). In~\cite{Po2}, Popa answered this question negatively.
 Indeed, Popa constructed
examples of maximal hyperfinite $II_1$ subfactors with relative
commutant isomorphism to $\mathbb{C}^n$ for any $n\geq 1$ and
hyperfinite $II_1$ subfactors with noncommutative relative
commutant. In~\cite{Ge}, Ge constructed  a maximal hyperfinite
$II_1$ subfactor of a $II_1$ factor with a non-injective relative
commutant!

The following lemmas show the relation between maximal injective von
Neumann subalgebras and maximal injective subfactors.
\begin{Lemma} If $\M$ is a factor and $\R$ is a maximal injective
subfactor of $\M$, then $\R$ is a maximal injective von Neumann
subalgebra of $\M$ if and only if $\R$ is irreducible in $\M$, i.e,
$\R^{'}\cap\M=\mathbb{C}I$.
\end{Lemma}
\begin{proof} If $\R$ is a maximal injective von Neumann subalgebra
of $\M$, then by Lemma 1.1, the center $\Z$ of $\R$ is
$\R^{'}\cap\M$. Since $\R$ is a factor, $\R^{'}\cap\M=\mathbb{C}I$.
Conversely, suppose $\R^{'}\cap\M=\mathbb{C}I$. For any injective
von Neumann algebra $\BB$ such that $\R\subseteq\BB\subseteq\M$, we
have $\BB^{'}\cap\BB\subseteq \R^{'}\cap\M=\mathbb{C}I$. Therefore,
$\BB$ is an injective subfactor of $\M$. Since $\R$ is a maximal
injective  subfactor of $\M$, $\R=\BB$.
\end{proof}

\begin{Lemma} Let $\M$ be a  factor and
$\BB$ be a maximal injective von Neumann subalgebra of $\M$. Let
$\Z$ be the center of $\BB$. If $\Z$ is atomic and $P_1,P_2,\cdots,$
are minimal  projections in $\Z$ such that $\sum P_i=I$. Then
$\BB_i=\BB_{P_i}$ is a maximal injective subfactor of $\M_{P_i}$
such that $\BB_i^{'}\cap\M_{P_i}=\mathbb{C}P_i$ for all $i$.
\end{Lemma}
\begin{proof} Since $\Z$ is atomic and $P_1,P_2,\cdots,$ are
minimal projections in $\Z$, $\BB_i$ is a subfactor of $\M_{P_i}$.
Since $\BB$ is injective, $\BB_i$ is injective. If $\L_i$ is an
injective von Neumann algebra such that $\BB_i\subseteq\L_i\subseteq
\M_{P_i}$, then $\BB_1\oplus\cdots\oplus\L_i\oplus\cdots$ is an
injective von Neumann subalgebra of $\M$ such that $\BB\subseteq
\BB_1\oplus\cdots\oplus\L_i\oplus\cdots\subseteq\M$. Since $\BB$ is
a maximal injective von Neumann subalgebra of $\M$,
$\BB=\BB_1\oplus\cdots\oplus\L_i\oplus\cdots$. This implies that
$\BB_i=\L_i$. So $\BB_i$ is a maximal injective von Neumann
subalgebra of $\M_{P_i}$. Since $\BB_i$ is a factor, $\BB_i$ is
irreducible in $\M_{P_i}$ by Lemma 1.5.
\end{proof}

\subsection{On tensor products of von Neumann algebras}
In~\cite{G-K}, Ge and Kadison proved the following basic theorem for
tensor products of von Neumann algebras.

\noindent{\bf Ge-Kadison's Splitting Theorem} \emph{If $\M_1$ is a
factor and $\M_2$ is a von Neumann algebra, and $\M$ is a von
Neumann subalgebra of $\M_1\otime\M_2$ which contains $\M_1\otime
\mathbb{C}I$, then $\M=\M_1\otime \N_2$ for some $\N_2$, a von
Neumann subalgebra of $\M_2$.}

Slice-map technique of Tomiyama~\cite{To2}  plays a key role in the
proof of Ge-Kadison's splitting theorem. Let $\M_1$ and $\M_2$ be
von Neumann algebras. With $\rho$ in $(\M_1)_{\#}$ (the predual of
$\M_1$), $\sigma$ in $(\M_2)_{\#}$ and $T\in\M_1\otime\M_2$, the
mapping $\sigma\rightarrow (\rho\otime\sigma)(T)$ is a bounded
linear functional on $(\M_2)_{\#}$, hence, an element
$\Psi_{\rho}(T)$ in $\M_2$. Symmetrically, we construct an operator
$\Phi_{\sigma}(T)$ in $\M_1$. The mappings $\Psi_{\rho}$ and
$\Phi_{\sigma}$ are referred to as \emph{slice mappings} (of
$\M_1\otime\M_2$ onto $\M_2$ and $\M_1$ corresponding to  $\rho$ and
$\sigma$, respectively). Tomiyama's Slice Mapping Theorem~\cite{To2}
says that if $\N_1$ and $\N_2$ are von Neumann subalgebras of $\M_1$
and $\M_2$, respectively, and $T\in \M_1\otime\M_2$, then
$T\in\N_1\otime\N_2$ if and only if $\Phi_{\sigma}(T)\in\N_1$ and
$\Psi_{\rho}(T)\in\N_2$ for each $\sigma\in (\M_2)_{\#}$ and
$\rho\in (\M_1)_{\#}$. For a generalization of Ge-Kadison's
splitting theorem, we refer to~\cite{S-Z}.

The following lemma is well-known. For the sake of completeness, we
include the proof here.
\begin{Lemma} If $\M_1,\N_1$ are von Neumann algebras acting on a Hilbert
space $\H$,  and $\M_2,\N_2$ are von Neumann algebras acting on a
Hilbert space $\K$, then $(\M_1\otime\M_2)\cap
(\N_1{'}\otime\N_2{'})\\
=(\M_1\cap\N_1{'})\otime (\M_2\cap\N_2{'}).$
\end{Lemma}
\begin{proof} It is obvious that $(\M_1\cap\N_1{'})\otime
(\M_2\cap\N_2{'})\subseteq (\M_1\otime\M_2)\cap
(\N_1{'}\otime\N_2{'})$. Conversely, if $T\in (\M_1\otime\M_2)\cap
(\N_1{'}\otime\N_2{'})$, then $\Phi_\sigma(T)\in\M_1$ and
$\Phi_\sigma(T)\in\N_1{'}$ for any $\sigma\in\B(\K)_{\#}$ by
Tomiyama's slice-mapping theorem. Thus
$\Phi_\sigma(T)\in\M_1\cap\N_1{'}$. Similarly, for any
$\rho\in\B(\H)_{\#}$, $\Psi_{\rho}(T)\in \M_2\cap\N_2{'}$. By
Tomiyama's slice-mapping theorem, $T\in (\M_1\cap\N_1{'})\otime
(\M_2\cap\N_2{'})$. Therefore, $(\M_1\otime\M_2)\cap
(\N_1{'}\otime\N_2{'})=(\M_1\cap\N_1{'})\otime (\M_2\cap\N_2{'}).$
\end{proof}

\section{Induced conditional expectations}
\begin{Lemma} Let $\M$ be a von Neumann algebra and $\L,\N$ be von
Neumann subalgebras of $\M$ such that $\L\subseteq\N$. If ${\bf E}$
is a conditional expectation from $\M$ onto $\N$, then ${\bf E}$
induces a conditional expectation from $\L{'}\cap\M$ onto
$\L{'}\cap\N$.
\end{Lemma}
\begin{proof} $\forall S\in\L{'}\cap\M$ and $T\in\L$, $ST=TS$. Apply
the conditional expectation ${\bf E}$ to both sides of $ST=TS$ and
note that $\L\subseteq\N$. We have ${\bf E}(S)T=T{\bf E}(S)$. Thus
${\bf E}(S)\in\L{'}\cap\N$. Since $\L{'}\cap\N\subseteq\L{'}\cap\M$,
${\bf E}$ is a conditional expectation from $\L{'}\cap\M$ onto
$\L{'}\cap\N$ when ${\bf E}$ is restricted on $\L{'}\cap\M$.
\end{proof}

The following is another key lemma.

\begin{Lemma} Let $\AA_i$  be a von Neumann algebra acting on a Hilbert
space $\H_i$, and $\N_i$ be a von Neumann subalgebra of $\AA_i$,
$i=1,2$. Let $\AA$ be a von Neumann algebra such that
$\N_1\otime\N_2\subseteq\AA\subseteq \AA_1\otime\AA_2$. If there is
a conditional expectation ${\bf E}$ from $\AA_1\otime\AA_2$ onto
$\AA$, then ${\bf E}$ induces a conditional expectation from
$(\N_1{'}\cap\AA_1)\otime \AA_2$ onto
$((\N_1{'}\cap\AA_1)\otime\AA_2)\cap\AA$.
\end{Lemma}
\begin{proof}  By Lemma 1.7, $(\AA_1\otime\AA_2)\cap (\N_1\otime
\mathbb{C}I)'=(\AA_1\otime\AA_2)\cap (\N_1{'}\otime
\B(\K))=(\N_1{'}\cap\AA_1)\otime\AA_2$. By Lemma 2.1, ${\bf E}$
induces a conditional expectation from $(\N_1{'}\cap\AA_1)\otime
\AA_2$ onto $((\N_1{'}\cap\AA_1)\otime\AA_2)\cap\AA$.
\end{proof}

\begin{Corollary} Assume the conditions of Lemma 2.2 and
 $\N_1{'}\cap\AA_1=\mathbb{C}I$. Let
$\L_2=\{T\in\AA_2:\, I\otimes T\in\AA\}$. Then ${\bf E}$ induces a
conditional expectation from $\AA_2$ onto $\L_2$.
\end{Corollary}

As an application of Lemma 2.2 and Corollary 2.3, we give a new
proof of Ge-Kadison's splitting theorem in the case when $\M_1$ and
$\M_2$ are finite. Let $\N$ be a von Neumann algebra such that
$\M_1\otime\mathbb{C}I\subseteq\N\subseteq\M_1\otime\M_2$. Then
there is a \emph{normal} conditional expectation ${\bf E}$ from
$\M_1\otime\M_2$ onto $\N$. By Corollary 2.3, ${\bf E}$ induces a
conditional expectation, denoted by ${\bf E}_2$, from $\M_2$ onto
$\N_2\triangleq\{T\in\M_2:\, I\otimes T\in\N\}$. Now for any
$S\in\M_1, T\in\M_2$, we have ${\bf E}(S\otimes T)=S\otimes {\bf
E}_2(T)\in \M_1\otime\N_2$. Since ${\bf E}$ is normal, $\N={\bf
E}(\M_1\otime\M_2)\subseteq \M_1\otime\N_2$. Since
$\N\supseteq\M_1\otime\N_2$, $\N=\M_1\otime\N_2$.\\

As another application of Lemma 2.2 and Corollary 2.3, we give a new
proof of Theorem 6.7 of~\cite{S-Z}.

\begin{Lemma} Let $\A$ be an abelian von Neumann
algebra, and $\AA_2$ be a minimal injective von Neumann algebra
extension of a von Neumann algebra $\N$. Suppose $\AA_2$ has
separable predual. If $\AA$ is an injective von Neumann algebra such
that $\A\otime\N\subseteq\AA\subseteq\A\otime \AA_2$, then
$\AA=\A\otime\AA_2$.
\end{Lemma}
\begin{proof} We can assume that $\A$ and $\AA_2$ are von Neumann
algebras acting on Hilbert spaces $\H$ and $\K$, respectively, in
standard form. Then $\A^{'}=\A$ and $\K$ is a separable Hilbert
space. By $\A\otime\N\subseteq\AA\subseteq\A\otime \AA_2$, we have
$\A\otime\AA_2^{'}\subseteq\AA^{'}\subseteq\A\otime \N{'}$. Note
that $\AA_2^{'}$ is a maximal injective von Neumann subalgebra of
$\N {'}$. By Lemma 6.6 of [S-Z], $\AA^{'}=\A\otime \N {'}$.
Therefore, $\AA=\A\otime\AA_2$.
\end{proof}

Lemma 2.4 is almost obvious in the case when $\A$ is atomic. If $\A$
is diffuse, it is natural to consider direct integrals. The proof of
Lemma 6.6 of~\cite{S-Z} is based on direct integrals. It would be
interesting if there is a ``global proof" of Lemma 2.4. Is Lemma 2.4
true without the assumption that $\AA_2$ has separable predual?

\begin{Theorem} Let $\M_1$ be an injective von Nuemann algebra
and $\M_2$ be a von Neumann algebra with separable predual. If
$\BB_2$ is a maximal injective von Neumann subalgebra of $\M_2$,
  then $\M_1\otime\BB_2$ is a maximal injective von Neumann
  subalgebra of $\M_1\otime \M_2$.
  \end{Theorem}
\begin{proof} We can assume that $\M_1$ and $\M_2$ are von Neumann
algebras acting on Hilbert spaces $\H$ and $\K$, respectively. Then
$\K$ is a separable Hilbert space.  Let $\A$ be the center of
$\M_1$. Suppose $\BB$ is an injective von Neumann algebra such that
$\M_1\otime\BB_2\subseteq\BB\subseteq\M_1\otime\M_2$. Then we have
$\M_1^{'}\otime\BB_2^{'}\supseteq\BB^{'}\supseteq\M_1^{'}\otime\M_2^{'}$.
Since $\BB^{'}$ is an injective von Neumann subalgebra of
$\M_1^{'}\otime\BB_2^{'}$, there is a conditional expectation ${\bf
E}$ from $\M_1^{'}\otime\BB_2^{'}$ onto $\BB^{'}$.   By Lemma 2.2,
${\bf E}$ induces a conditional expectation from $\A\otime\BB_2^{'}$
onto $\AA\triangleq(\A\otime\BB_2^{'})\cap\BB^{'}$. So $\AA$ is an
injective von Neumann algebra such that
$\A\otime\BB_2^{'}\supseteq\AA\supseteq \A\otime\M_2^{'}$. Since
$\BB_2$ is a maximal injective von Neumann subalgebra of $\M_2$,
$\BB_2^{'}$ is a minimal injective von Neumann algebra extension of
$\M_2^{'}$. By Lemma 2.4, $\AA=\A\otime\BB_2^{'}$. Thus
$\mathbb{C}I\otime \BB_2^{'}\subseteq\AA\subseteq\BB^{'}$. So
$\BB^{'}=\M_1^{'}\otime\BB_2^{'}$
 and
$\BB=\M_1\otime\BB_2$.
 \end{proof}

 \section{Popa's question in the case when $\BB_1$ is a factor}

 \begin{Theorem} Let $\M_i$ be a von Neumann algebra,
  and $\BB_i$ be a maximal
injective von Neumann subalgebra of $\M_i$, $i=1,2$. If $\M_1$ has
separable predual and $\BB_1$ is a factor, then $\BB_1\otime\BB_2$
is a maximal injective von Neumann subalgebra of $\M_1\otime \M_2.$
\end{Theorem}
\begin{proof}
 Let $\BB$ be an injective von Neumann algebra
such that $\BB_1\otime\BB_2\subseteq\BB\subseteq \M_1\otime\M_2$. To
prove the theorem, we need to show that $\BB=\BB_1\otime\BB_2$. We
can assume that $\M_1$ and $\M_2$ are von Neumann algebras acting on
Hilbert spaces $\H$ and $\K$, respectively. So we have
$\BB_1^{'}\otime\BB_2^{'}\supseteq\BB^{'}\supseteq
\M_1^{'}\otime\M_2^{'}$.    Since $\BB_1,\BB_2,\BB$ are injective,
$\BB_1^{'},\BB_2^{'},\BB^{'}$ are injective. Since $\BB_i$ is a
maximal injective von Neumann algebra of $\M_i$, $\BB_i^{'}$ is a
minimal injective von Neumann algebra extension of $\M_i^{'}$,
$i=1,2.$\\

Since $\BB^{'}$ is an injective von Neumann subalgebra of
$\BB_1^{'}\otime\BB_2^{'}$, there is a conditional expectation ${\bf
E}$ from $\BB_1^{'}\otime\BB_2^{'}$ onto $\BB^{'}$. Let
$\L_2=\{T\in\BB_2^{'}:\, 1\otimes T\in \BB'\}$. Then
$\L_2\subseteq\BB_2^{'}$. By Lemma 1.1,
$(\M_1^{'})^{'}\cap\BB_1^{'}=\BB_1^{'}\cap\M_1=\BB_1^{'}\cap
\BB_1=\mathbb{C}I$. By Corollary 2.3, ${\bf E}$ induces a
conditional expectation from $\BB_2^{'} $ onto $\L_2 $. Thus $\L_2$
is injective. Since $\M_2^{'}\subseteq \L_2\subseteq \BB_2^{'}$ and
$\BB_2^{'}$ is a minimal injective von Neumann algebra extension of
$\M_2^{'}$, $\L_2=\BB_2^{'}$. So $\mathbb{C}I\otime
\BB_2^{'}\subseteq \BB'$. This implies that
$\BB_1^{'}\otime\BB_2^{'}\supseteq\BB^{'}\supseteq
\M_1^{'}\otime\BB_2^{'}$ and hence
$\BB_1\otime\BB_2\subseteq\BB\subseteq \M_1\otime\BB_2$.
 By Theorem 2.5, $\BB=\BB_1\otime\BB_2$.
\end{proof}

\section{Popa's question in the case when the center of  $\BB_1$ is atomic}
The following is the main result of this paper.
\begin{Theorem} Let $\M_i$ be a
von Neumann algebra, and $\BB_i$ be a maximal injective von Neumann
subalgebra of $\M_i$, $i=1,2$. If $\M_1$ has separable predual and
the center of $\BB_1$ is atomic, then $\BB_1\otime\BB_2$ is maximal
injective in $\M_1\otime\M_2$.
\end{Theorem}
\begin{proof}
 Let $\BB$ be an injective von Neumann algebra
such that $\BB_1\otime\BB_2\subseteq\BB\subseteq \M_1\otime\M_2$. We
can assume that $\M_1$ and $\M_2$ are von Neumann algebras acting on
Hilbert spaces $\H$ and $\K$, respectively. To prove the theorem, we
need to show that $\BB=\BB_1\otime\BB_2$.
  Using Theorem~2.5, it is sufficient to
prove that $\BB \subseteq \M_1 \otime \BB_2$ and so, it is
sufficient
to prove that $\mathbb{C}I \otime \BB_2' \subseteq \BB'$.\\

Denote $\A=\BB_1 \cap \BB_1' = \M_1 \cap \BB_1'$, which is atomic,
with the set of minimal projections $\{P_n:\,\,n=1,2,\cdots\}$. Set
$\AA= \BB' \cap (\A \otime \BB_2') = \BB' \cap (\M_1' \otime
\mathbb{C}I)'$. As in Section~2 of the paper, $\AA$ is injective and
$$\mathbb{C}I \otime \M_2' \subseteq \AA \subseteq \A \otime\BB_2' \; .$$
Note that $P_n \otimes I \in \AA'$ for every $n$. Since $\BB_2'$ is
a minimal injective von Neumann algebra extension of $\M_2'$, we
have $\AA(P_n \otimes I) = P_n \otimes \BB_2'$ for every $n$. Denote
by $Z_n$ the smallest projection in $\mathcal{Z}(\AA)$ satisfying
$P_n \otimes I \leq Z_n$. We get $^*$-isomorphisms $\theta_n :
\BB_2' \rightarrow \AA Z_n$ uniquely determined by the formula
$$\theta_n(Y) (P_n \otimes I) = P_n \otimes Y \quad\text{for all}\;\; Y \in
\BB_2' \; .$$ Since $\mathbb{C}I \otime \M_2' \subseteq \AA$, it
follows that $\theta_n(X) = (I \otimes X)Z_n$ for all $X \in
\M_2'$.\\

So, for every $n,m$ and all $X \in \M_2'$, we have $\theta_n(X)Z_m =
\theta_m(X)Z_n$. Since $\BB_2'$ is a minimal injective von Neumann
algebra extension of $\M_2'$, by Corollary 1.4, the same formula
holds for all $Y \in \BB_2'$ and all $n,m$. This compatibility
formula yields for every $Y \in \BB_2'$ an element $A \in \AA$ such
that $AZ_n = \theta_n(Y)$ for all $n$. In particular, $A(P_n \otimes
I) = P_n \otimes Y$, i.e., \ $A=I \otimes Y$. So, $\mathbb{C}I
\otime \BB_2' \subseteq \AA$, ending the proof.
\end{proof}

 Replacing the use of the minimal
projections $P_n$ by a careful
  analysis of \lq infinitesimal projections\rq\ (i.e.\ using direct
  integral techniques), the same kind of idea maybe allows to prove
  the general case, not assuming $\mathcal{Z}(\BB_1)$ to be atomic.

\section{A result on maximal injective subfactor of tensor products of von Neumann algebras}
\begin{Theorem} Let $\M_i$ be a factor,
 and $\R_i$ be a maximal injective subfactor of $\M_i$, $i=1,2$.
 If $\R_1^{'}\cap\M_1\simeq\mathbb{C}^N$ \emph{(}$1\leq N\leq \infty$\emph{)} and
$\R_2^{'}\cap\M_2=\mathbb{C}I$, then $\R_1\otime\R_2$ is a maximal
injective subfactor of $\M_1\otime \M_2.$
\end{Theorem}
\begin{proof}
Consider an injective factor $\R$ such that
$\R_1\otime\R_2\subseteq\R\subseteq \M_1\otime\M_2$. We need to show
that $\R=\R_1\otime\R_2$. We can assume that $\M_1$ and $\M_2$ are
von Neumann algebras acting on Hilbert spaces $\H$ and $\K$,
respectively. So we have
$\R_1^{'}\otime\R_2^{'}\supseteq\R^{'}\supseteq
\M_1^{'}\otime\M_2^{'}$. By Lemma 1.5, $\R_2$ is a maximal injective
von Neumann subalgebra of $\M_2$ and thus $\R_2^{'}$ is a minimal
injective von Neumann algebra extension of $\M_2^{'}$. By
assumption, $\R_2^{'}\cap\M_2=\R_2^{'}\cap\R_2=\mathbb{C}I$.\\

Let ${\bf E}$ be a conditional expectation from
$\R_1^{'}\otime\R_2^{'}$ onto $\R^{'}$, and
$\A=\R_1^{'}\cap\M_1\simeq\mathbb{C}^N$. By Lemma 2.2, ${\bf E}$
induces a conditional expectation from $\A\otime \R_2^{'}$ onto
$\BB\triangleq(\A\otime\R_2^{'})\cap\R^{'}$. Therefore, $\BB$ is an
injective von Neumann algebra such that
$\R'\supseteq\BB\supseteq(\A\otime\R_2^{'})\cap(\M_1^{'}\otime\M_2^{'})=
\mathbb{C}I\otime\M_2^{'}$.\\

Similar arguments as the proof of Theorem 4.1 show that
$\BB\supseteq\mathbb{C}I\otime \R_2^{'}$. So $\R'\supseteq
\mathbb{C}I\otime \R_2^{'}$.  By Ge-Kadison's splitting theorem (see
1.4), $\R^{'}=\N^{'}\otime
 \R_2^{'}$ for some von Neumann subalgebra $\N^{'}$ of $\R_1^{'}$.
 Therefore $ \R=\N\otime\R_2$.
 Since $\R$ is an injective factor,
$\N$ is an injective factor such that
$\R_1\subseteq\N\subseteq\M_1$. Since $\R_1$ is a maximal injective
subfactor of $\M_1$, $\N=\R_1$.  Therefore, $\R=\R_1\otime\R_2$.\\

\end{proof}

\begin{Corollary}Let $\M_i$ be a factor and $\R_i$ be an
irreducible, maximal injective subfactor of $\M_i$, $i=1,2$. Then
$\R_1\otime\R_2$ is an irreducible, maximal injective subfactor of
$\M_1\otime \M_2$.
\end{Corollary}

\section{Concluding remarks}
\subsection{}  Let $\M_i$ be a von Neumann algebra and
$\BB_i$ be a maximal injective von Neumann subalgebra of $\M_i$,
$i=1,2$.  Suppose $\M_1$ is a type $II_1$ von Neumann algebra with
separable predual and $\BB_1$ is an maximal injective von Neumann
subalgebra of $\M_1$.  By~\cite{Co}, $\BB_1=(\A\otime \R)\bigoplus
\oplus_{n=1}^{\infty}(\A_n\otime M_n(\mathbb{C}))$, where $\A,
\A_1,\A_2, \cdots, $ are abelian von Neumann algebras and $\R$ is
the hyperfinite type $II_1$ factor. By Theorem 4.1, if
$\BB_1=\A\otime \R$ is type $II_1$ and $\A$ is atomic, then
$\BB_1\otime\BB_2$ is a maximal injective von Neumann subalgebra of
$\M_1\otime\M_2$. Popa's question remains open for all other cases
of $\BB_1$, e.g., $\BB_1$ is abelian. It is still not known that  a
diffuse abelian von Neumann algebra with separable predual can be
embedded into any  non hyperfinite separable type $II_1$ factor as a
maximal injective von Neumann subalgebra or not. Recently,
J.~Shen~\cite{Sh} proved that $\{\L(a)\}^{''}\otime\{\L(a)\}^{''}$
is a maximal injective von Neumann algebra of
$\L(\mathbb{F}_n)\otime\L(\mathbb{F}_n)$. J.~Shen also provided the
first example of a  Mcduff $II_1$ factor which contains an abelian
von Neumann algebra as a maximal injective von Neumann subalgebra.
\subsection{}  Let $\M_i$ be a von Neumann algebra acting on a Hilbert
space $\H_i$, and
 $\BB_i$ be a maximal injective von Neumann subalgebra
of  $\M_i$, $i=1,2$. Suppose $\BB$ is an injective von Neumann
algebra such that
$\BB_1\otime\BB_2\subseteq\BB\subseteq\M_1\otime\M_2$.  Then we have
$\BB_1^{'}\otime\BB_2^{'}\supseteq\BB^{'}\supseteq\M_1^{'}\otime\M_2^{'}$.
Let ${\bf E}$ be a conditional expectation from
$\BB_1^{'}\otime\BB_2^{'}$ onto $\BB^{'}$. By Lemma 2.2, ${\bf E}$
induces a conditional expectation from $\Z\otime \BB_2^{'}$ onto
$\AA\triangleq(\Z\otime \BB_2^{'})\cap\BB^{'}$, where
$\Z=\BB_1^{'}\cap\BB_1=\BB_1^{'}\cap \M_1$. Note that
$\Z\otime\BB_2^{'}\supseteq\AA\supseteq \mathbb{C}I\otime
\M_2^{'}$. This leads to the following question:\\

{\bf Question $1$:} Suppose $\A$ is an abelian von Neumann algebra
and $\AA_2$ is a minimal injective von Neumann algebra extension of
a von Neumann algebra $\N_2$. If $\AA$ is an injective von Neumann
algebra such that
$\A\otime\AA_2\supseteq\AA\supseteq\mathbb{C}I\otime \N_2$, is
$\AA\supseteq\mathbb{C}I\otime \AA_2$?\\

An affirmative answer to Question 1 would give rise to an
affirmative answer to Popa's question (with assumption that $\M_1$
has separable predual). Indeed, if $\AA\supseteq \mathbb{C}I\otime
\BB_2^{'}$, then $\BB^{'}\supseteq \AA\supseteq \mathbb{C}I\otime
\BB_2^{'}$. Therefore,
$\BB_1^{'}\otime\BB_2^{'}\supseteq\BB^{'}\supseteq\M_1^{'}\otime\BB_2^{'}$.
Hence, $\BB_1\otime\BB_2\subseteq \BB\subseteq \M_1\otime \BB_2$.
Apply Theorem 2.5, $\BB=\BB_1\otime\BB_2$. So $\BB_1\otime\BB_2$ is
a maximal injective von Neumann subalgebra
of $\M_1\otime\M_2$.\\

In Question 1, we may  assume that $\A$ and $\AA_2$ are von Neumann
algebras acting on Hilbert spaces $\H$ and $\K$, respectively, in
standard form. Then $\A$ is a maximal abelian von Neumann subalgebra
of $\B(\H)$, i.e., $\A^{'}=\A$.
 Consider the  commutant of $\A\otime\AA_2,$ $\AA,$ $\mathbb{C}I\otime
\N_2$, respectively, Question 1 is equivalent to the following
question:\\

{\bf Question $1'$:} Suppose $\A$ is a maximal abelian von Neumann
algebra and $\BB_2$ is a maximal injective von Neumann subalgebra of
a von Neumann algebra $\M_2$.  If $\BB$ is an injective von Neumann
algebra such that
$\A\otime\BB_2\subseteq\BB\subseteq\B(\H)\otime\M_2$, is
$\BB\subseteq\B(\H)\otime\BB_2$?\\

By the proof of Theorem 4.1, if $\A$ is atomic, then the answer to
question 1 is affirmative and thus to question $1'$ is affirmative.
\vskip 1cm

{\bf Acknowledgements.}\,\, The author thanks the referee for a
simple proof of Theorem 4.1 and other useful remarks that greatly
improve the exposition of the paper.

\end{document}